\documentclass[twoside]{article}

\oddsidemargin=\evensidemargin
\catcode`\@=11
\long\def\@makefntext#1{
\protect\noindent \hbox to 3.2pt {\hskip-.9pt  
$^{{\eightrm\@thefnmark}}$\hfil}#1\hfill}		%CAN BE USED 

\def\ps@myheadings{\let\@mkboth\@gobbletwo		%SIZE OF R/H NOS.
\def\@oddhead{\hbox{}
\rightmark\hfil\eightrm\thepage}   
\def\@oddfoot{}\def\@evenhead{\eightrm\thepage\hfil
\leftmark\hbox{}}\def\@evenfoot{}
\def\sectionmark##1{}\def\subsectionmark##1{}}

\catcode`@=11		      		     %SIZE OF OPENING PAGE NO.
\def\ps@plain{\let\@mkboth\@gobbletwo
     \def\@oddhead{}\def\@oddfoot{\eightrm\hfil\thepage
     \hfil}\def\@evenhead{}\let\@evenfoot\@oddfoot}

\newcounter{sectionc}\newcounter{subsectionc}\newcounter{subsubsectionc}
\renewcommand{\section}[1] {\vspace{12pt}\addtocounter{sectionc}{1} 
\setcounter{subsectionc}{0}\setcounter{subsubsectionc}{0}\noindent 
	{\tenbf\thesectionc. #1}\par\vspace{5pt}}
\renewcommand{\subsection}[1] 
{\vspace{12pt}\addtocounter{subsectionc}{1} 
	\setcounter{subsubsectionc}{0}\noindent 
	{\bf\thesectionc.\thesubsectionc. 
	{\kern1pt \bfit #1}}\par\vspace{5pt}}
\renewcommand{\subsubsection}[1] {\vspace{12pt}
	\addtocounter{subsubsectionc}{1}
	\noindent
	{\tenrm\thesectionc.\thesubsectionc.\thesubsubsectionc.	{\kern1pt 
	\it #1}}\par\vspace{5pt}}
\newcommand{\nonumsection}[1] {\vspace{12pt}\noindent{\tenbf #1}
	\par\vspace{5pt}}

\newcommand{\textlineskip}{\baselineskip=13pt}
\newcommand{\smalllineskip}{\baselineskip=10pt}

\newcommand{\copyrightheading}[1]
	{\vspace*{-2.5cm}\smalllineskip{\flushleft
	{\footnotesize Journal of Knot Theory and Its Ramifications #1}\\
   	{\footnotesize \copyright\kern2pt World Scientific 
         Publishing Company}\\
         }}

\def\abstracts#1#2#3#4{{
	\centering{\begin{minipage}{4.5in}\footnotesize\baselineskip=10pt
	\centerline{ABSTRACT} 
	\parindent=15pt #1\par 
	\parindent=15pt #2\par
	\parindent=15pt #3\par
	\parindent=15pt #4\par
	\end{minipage}}\par}} 

\def\keywords#1{{ 
	\centering{\begin{minipage}{4.5in}\footnotesize\baselineskip=10pt
	{\footnotesize\it Keywords}\/: #1
	\end{minipage}}\par}}

\renewenvironment{thebibliography}[1]
	{\frenchspacing

	 \ninerm\baselineskip=11pt
	 \begin{list}{[\arabic{enumi}]}
	{\usecounter{enumi}\setlength{\parsep}{0pt}
	 \setlength{\leftmargin 13.7pt}{\rightmargin 0pt} %[1--9] ITEMS
	 \setlength{\itemsep}{0pt} \settowidth
	{\labelwidth}{[#1]}\sloppy}}{\end{list}}

\newcounter{itemlistc}
\newcounter{romanlistc}
\newcounter{alphlistc}
\newcounter{arabiclistc}

\newenvironment{alphlist}
	{\setcounter{alphlistc}{0}
	 \begin{list}{$($\alph{alphlistc}$)$}
	{\usecounter{alphlistc}
	 \setlength{\parsep}{0pt}
	 \setlength{\itemsep}{0pt}}}{\end{list}}

\newcommand{\fcaption}[1]{
        \refstepcounter{figure}
        \setbox\@tempboxa = \hbox{\footnotesize Fig.~\thefigure. #1}
        \ifdim \wd\@tempboxa > 5in
           {\begin{center}
        \parbox{5in}{\footnotesize\smalllineskip Fig.~\thefigure. #1}
            \end{center}}
        \else
             {\begin{center}
             {\footnotesize Fig.~\thefigure. #1}
              \end{center}}
        \fi}

\def\pmb#1{\setbox0=\hbox{#1}
	\kern-.025em\copy0\kern-\wd0
	\kern.05em\copy0\kern-\wd0
	\kern-.025em\raise.0433em\box0}

\def\fnt#1#2{\footnotetext{\kern-.3em
	{$^{\mbox{\scriptsize #1}}$}{#2}}}

\def\fpage#1{\begingroup
\voffset=.3in
\thispagestyle{empty}\begin{table}[b]\centerline{\footnotesize #1}
	\end{table}\endgroup}

\def\runninghead#1#2{\pagestyle{myheadings}
\markboth{{\protect\footnotesize\it{\quad #1}}\hfill}
{\hfill{\protect\footnotesize\it{#2\quad}}}}

\font\tenrm=cmr10
\font\tenit=cmti10 
\font\tenbf=cmbx10
\font\bfit=cmbxti10 at 10pt
\font\ninerm=cmr9
\font\nineit=cmti9
\font\ninebf=cmbx9
\font\eightrm=cmr8

\def\@begintheorem#1#2{\trivlist	%6/9/94
	\item[\hskip\labelsep{\bf #1\ #2.}]} 
\def\@opargbegintheorem#1#2#3{\trivlist
	\item[\hskip\labelsep{\bf #1\ #2\ (#3).}]}

\newenvironment{itemlist1}			%WILL GENERATE 
    	{\setcounter{itemlistc}{0}		%\NOINDENT\BULLET 
	 \begin{list}{$\bullet$}		%AUTOMATICALLY
	{\usecounter{itemlistc}			%UNLESS OTHERWISE SPECIFIED
	 \leftmargin10pt	       %IE.  0PT = MOVE OUTSIDE LEFT MARGIN
	 \setlength{\parsep}{0pt}
	 \setlength{\itemsep}{0pt}     %ADD EXTRA VSPACE BETWEEN ITEMS
	}}{\end{list}}
\def\qed{\hbox{${\vcenter{\vbox{			%HOLLOW SQUARE
   \hrule height 0.4pt\hbox{\vrule width 0.4pt height 6pt
   \kern5pt\vrule width 0.4pt}\hrule height 0.4pt}}}$}}

\newcommand{\R}{\rm I\kern-.2em R}
\usepackage{epsfig}
\begin{document}
\setlength{\textheight}{7.7truein}  %for 2nd page onwards

\runninghead{Computing the Writhe of a Knot}
{Computing the Writhe of a Knot}

\normalsize\textlineskip
\thispagestyle{empty}
\setcounter{page}{1}

\copyrightheading{}		    %{Vol.~0, No.~0 (1999) 00--00}

\vspace*{0.88truein}

\fpage{1}
\centerline{\bf COMPUTING THE WRITHE OF A KNOT}
\baselineskip=13pt
\vspace*{0.37truein}
\centerline{\footnotesize DAVID CIMASONI}
\baselineskip=12pt
\centerline{\footnotesize\it Section de math\'ematiques}
\baselineskip=10pt
\centerline{\footnotesize\it Universit\'e de Gen\`eve}
\baselineskip=10pt
\centerline{\footnotesize\it e-mail: David.Cimasoni@math.unige.ch}

\vspace*{1truein}

\abstracts{We study the variation of the Tait number of a closed space curve
according to its different projections. The results are used to compute the writhe of a knot,
leading to a closed formula in case of polygonal curves.}{}{}{}

\vspace*{10pt}

\vspace*{1pt}\textlineskip
\section{Introduction}
\vspace*{-0.5pt}
Since Crick and Watson's celebrated article in 1953, the local structure of DNA is well 
understood; its visualization as a double helix suggests the mathematical model of a ribbon 
in $\R^3$. But what about the global structure of DNA, that is, what kind of closed curve does 
the core of the ribbon form? It appears that these curves can show a great complexity 
(supercoiling) and are of central importance in understanding DNA (see [1] and [2]).

In 1961, C\u alug\u areanu [3] made the following discovery: take a ribbon in ${\R}^3$, let $Lk$ be the 
linking number of its border components, and $Tw$ its total twist; then the difference 
$Lk - Tw$  depends only on the core of the ribbon. This real number, later called 
{\it writhe} by Fuller [4], is of great interest for biologists, as it gives a measure of supercoiling in 
DNA.

Several techniques have been developed to estimate the writhe of a given space curve, 
particularly by Aldinger, Klapper and Tabor [5]. The aim of the present paper is to give a 
new way of computing the writhe (see proposition 3). Moreover, this method will lead to a closed formula for 
the writhe of any polygonal space curve.

\section{Definitions}
\noindent{\bf Definition~1.} {\it A polygonal knot (or PL knot) is the union of a finite 
number of segments in ${\R}^3$, homeomorphic to $S^1$. A point $x \in K$ is either a vertex 
or an interior point of $K$.}\medskip

Let us fix a PL knot $K$. Take a point $\xi$ in $S^2$; let $d_{\xi}$ be the oriented 
vector line containing $\xi$ and $p_{\xi}\colon K \to {\R}^2$ the orthogonal projection with 
kernel $d_{\xi}$.\medskip

\noindent {\bf Definition~2.} {\it The map $p_{\xi}$ is a good projection (or generic 
projection) if, for all $v$ in ${\R}^2$, $p_\xi^{-1}(v)$ is empty, one point, or two interior 
points of $K$.}\medskip

Clearly, the good projections form an open dense subset ${\cal O}_K$ of $S^2$; given $\xi \in 
{\cal O}_K$, we get a diagram of the PL knot $K$ by pointing out which segment lies on the 
top of the other at double points, according to the orientation of $d_{\xi}$.

\begin{figure}[htbp] 
\vspace*{13pt}
\centerline{\psfig{file=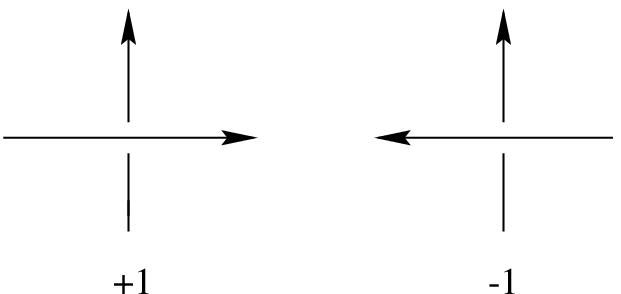,height=2.5cm}}
\vspace*{13pt}
\fcaption{}
\label{fig:Tait}
\end{figure}

To each double point $v$ of a diagram, it is possible to give a sign $s(v)=\pm 1$ by the 
following rule: we (temporarily) orientate $K$, and set $s(v)=+1$ if the oriented upper 
strand has to be turned counterclockwise to coincide with the lower one, $s(v)=-1$ in the 
other case (see figure \ref{fig:Tait}). The sign obviously does not depend on the 
orientation of the knot.\medskip

\noindent {\bf Definition~3.} {\it Let $K$ be a PL knot and $\xi \in {\cal O}_K$; the Tait
number of K relatively to $\xi$ is the integer $T_K(\xi)= \sum_v s(v)$, where the sum runs 
through all double points of the diagram associated with $\xi$.}\medskip

Given $K$, let us define a function $T_K \colon S^2 \to {\R}$ by
$$T_K(\xi)=\cases{ T_K(\xi), & if $\xi \in {\cal O}_K$; \cr           
            	 0, & otherwise.\cr}$$
This function is not continuous, but it is integrable (see proposition 1).\medskip

\noindent {\bf Definition~4.} {\it The writhing number (or writhe) of a PL knot K is the real 
number $Wr(K) = \frac{1}{4\pi} \int\!\!\!\int_{S^2}T_K(\xi)\,d\xi$.}\medskip
         	 
Since $S^2 \backslash {\cal O}_K$ is of measure zero, we can extend the Tait number in 
any way without changing the value of the writhe; the first proposition will give a very 
natural way to do so.\medskip

The study of $T_K(\xi)$ requires one more object: the indicatrix of a knot, that we define now.
Let $K$ be a $PL$ knot; an orientation of $K$ allows us to number its segments $S_1, \dots,S_n$. 
For all $i\/$, $S_i$ determines a direction, that is, a point $s_i$ in $S^2$. For $1\le i \le
n-1$, let $\Gamma_i$ be the oriented arc of great circle on $S^2$ joining $s_i$ to $s_{i+1}$, 
and $\Gamma_n$ the one joining  $s_n$ to $s_1$. The union $\Gamma=\cup_{i=1}^n \Gamma_i$ is a
closed oriented curve on $S^2$. Of course, the opposite orientation of $K$ produces the antipodal
curve $-\Gamma$.\medskip

\noindent {\bf Definition~5.} {\it The (spherical) indicatrix $I$ of a $PL$ knot is the 
oriented curve of $S^2$ given by $I=\Gamma \cup -\Gamma$.}\medskip

It is important to notice that the indicatrix may run several times on a fixed arc of great 
circle. 

We can easily check that $p_\xi$ is locally injective if and only if $\xi \in S^2 \backslash I$. 

\section{The Results}
\noindent {\bf Proposition~1.} {\it If $\xi$ and $\xi \prime$ in ${\cal O}_K$ belong to the same 
connected component of $S^2\backslash I$, then $T_K(\xi)=T_K(\xi \prime)$.}\medskip

\noindent {\bf Remarks.}
\begin{itemlist1}
\item It is possible to extend the definition of the Tait number to all
locally injective projections in a very natural way by setting $T_K(\xi) = T_K(A)$, where 
$T_K(A)$ is the value of the Tait number on the connected component $A$ of $S^2\backslash I$ 
containing $\xi$.
\item This proposition directly implies that $T_K$ is integrable on $S^2$.
\end{itemlist1}

\noindent{\bf Proof.} Let $\Omega$ be a connected component of $S^2\backslash I$, 
and $\xi,\xi\prime \in \Omega \cap {\cal O}_K$. If $\Omega \cap {\cal O}_K$ were connected, 
a simple continuity argument would do, but this is not the case in general. Indeed, 
$\Omega \backslash (\Omega \cap {\cal O}_K)$ is the union of a finite number of curves on 
$S^2$ that can separate $\xi$ and $\xi\prime$. We need to describe all these curves, and 
check that crossing them does not change the value of $T_K$.

\begin{figure}[htbp]
\vspace*{13pt}
\centerline{\psfig{file=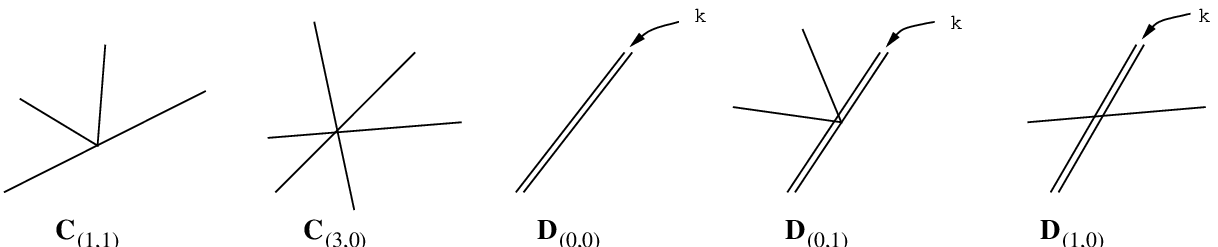}} 
\vspace*{13pt}
\fcaption{}
\label{fig:cinq}
\end{figure}

Let us first consider the non-generic projections that do not send coplanar segments onto the 
same line; we call $C_{(m,n)}$ the set of points $\xi \in S^2$ such that there exists $v \in 
{\R}^2$ with $p_\xi^{-1}(v)$ consisting of $m$ interior points and $n$ vertices of $K$. It is 
easy to check that if $m \ge 4$ or $n \ge 2$, $C_{(m,n)}$ is a finite number of points, as 
well as  $C_{(2,1)}$ and $C_{(3,1)}$. Since $C_{(0,0)}$, $C_{(1,0)}$ and $C_{(2,0)}$ are the
constituents of ${\cal O}_K$, the only potential trouble makers are $C_{(1,1)}$ and 
$C_{(3,0)}$ (see figure \ref{fig:cinq}).

We also have to study projections sending $k \ge 2$ coplanar (pairwise non-adjacent)
segments of $K$ onto a line $d$ of the plane. Let $D_{(m,n)}$ be the set of $\xi \in S^2$ such 
that there exists $v \in d$ with $p_\xi^{-1}(v)$ consisting of $m$ interior points and $n$ 
vertices of $K$. This time, only $D_{(0,0)}$, $D_{(0,1)}$ and $D_{(1,0)}$ need to be considered.

Thus, there are five types of non-generic projections that can separate $\xi$ and $\xi\prime$;
they are illustrated in figure \ref{fig:cinq}.

We now study $\Delta T = T_K(\eta)-T_K(\eta\prime)$, where $\eta$ and $\eta\prime$ are 
good projections on each side of a curve in $S^2$ formed by one type of non-generic projection.
Crossing $C_{(1,1)}$ (resp. $C_{(3,0)}$) corresponds to the second (resp. third) Reidemeister 
move; the Tait number being unchanged by these transformations, the first two cases are 
settled. It is trivial that $\Delta T=0$ for $D_{(1,0)}$, while $D_{(0,1)}$ can be seen as 
$k$ second Reidemeister moves.

\begin{figure}[htbp] 
\vspace*{13pt}
\centerline{\psfig{file=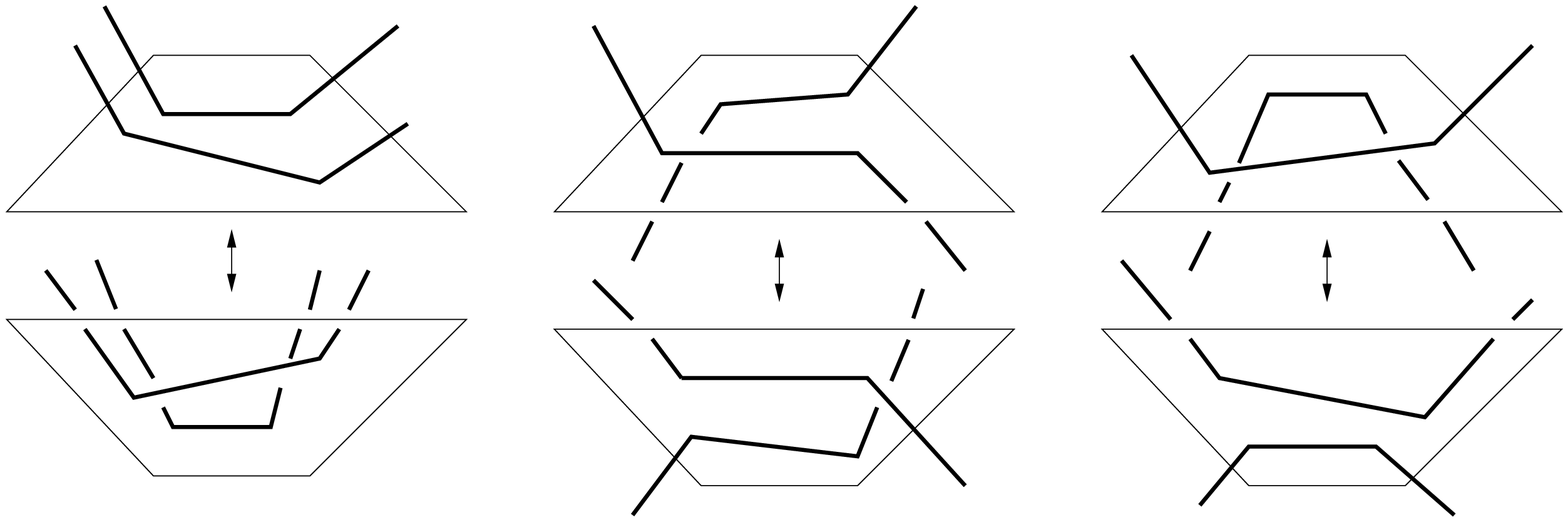, width=10cm}}
\vspace*{13pt}
\fcaption{}
\label{fig:D}
\end{figure}

It remains to show that $\Delta T=0$ for $D_{(0,0)}$. If $k=2$, all the possible cases are 
described in figure \ref{fig:D}; each time, $\Delta T=0$. For $k \ge 3$, the segments being pairwise
non-adjacent, we can apply the same argument (k-1) times. This concludes the proof.
\hfill \qed\kern0.8pt \medskip

What about the behavior of $T_K(\xi)$ when $\xi$ crosses the indicatrix? Let $\alpha$ be an 
open segment of the indicatrix $I$ on $S^2$, $p_1$ and $p_2$ two distinct points on $\alpha$. 
Let $n$ be the algebraic number of times that the indicatrix runs from $p_1$ to $p_2$. We 
will say that $\xi \in S^2\backslash \alpha$ is to the north of $\alpha$ if $(p_1 \times p_2) \cdot \xi$ 
is strictly positive, to the south if it is strictly negative.\medskip

\noindent {\bf Proposition~2.} {\it Let $\Omega_0$ be the component of 
$S^2\backslash I$ to the south of $\alpha$, $\Omega_1$ to the north; then: 
$T_K(\Omega_1)=T_K(\Omega_0)+n$.}\medskip

\noindent{\bf Proof.} To simplify the exposition, we will give the demonstration only when 
$\alpha$ is covered once by the indicatrix. 

\begin{figure}[htbp]
\vspace*{13pt}
\centerline{\psfig{file=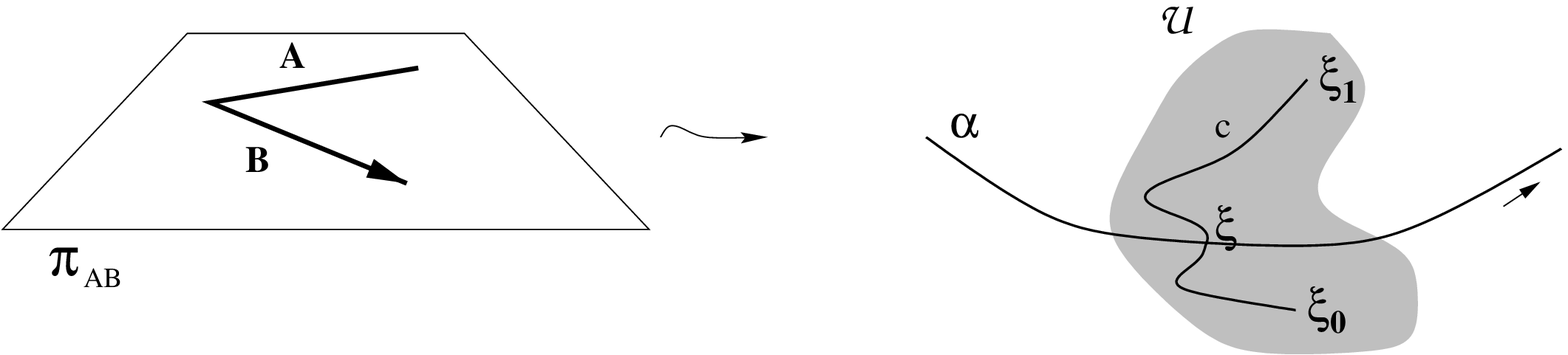,width=12cm}}
\vspace*{13pt}
\fcaption{}
\label{fig:prop2}
\end{figure}

Since $K$ is polygonal, $I$ is the union of a finite number of arcs of great circles; $\alpha$ 
is an open arc of great circle produced by two adjacent segments $A$ and $B$ of $K$ (we will 
say: the site $AB$). Let $U$ be an open path-connected subset of $S^2$ such that $U \backslash 
(U \cap \alpha) \subset {\cal O}_K$, and let us take $\xi_0 \in \Omega_0 \cap U, \xi_1 \in 
\Omega_1 \cap U$, $c$ a path in $U$ joining $\xi_0$ to $\xi_1$, crossing $\alpha$ at $\xi$ 
(see figure \ref{fig:prop2}). 
Since the indicatrix runs only once through $\alpha$ and since $Im(c) \backslash \xi \subset
{\cal O}_K$, the site $AB$ alone influences the variation $\Delta T$ of the Tait number 
along $c$.

Thus, we have to look at the site projected by $p_{\xi_0}$ and by $p_{\xi_1}$, and check that 
$T(\xi_1)=T(\xi_0)+1$. 

\begin{figure}[htbp]
\vspace*{13pt}
\centerline{\psfig{file=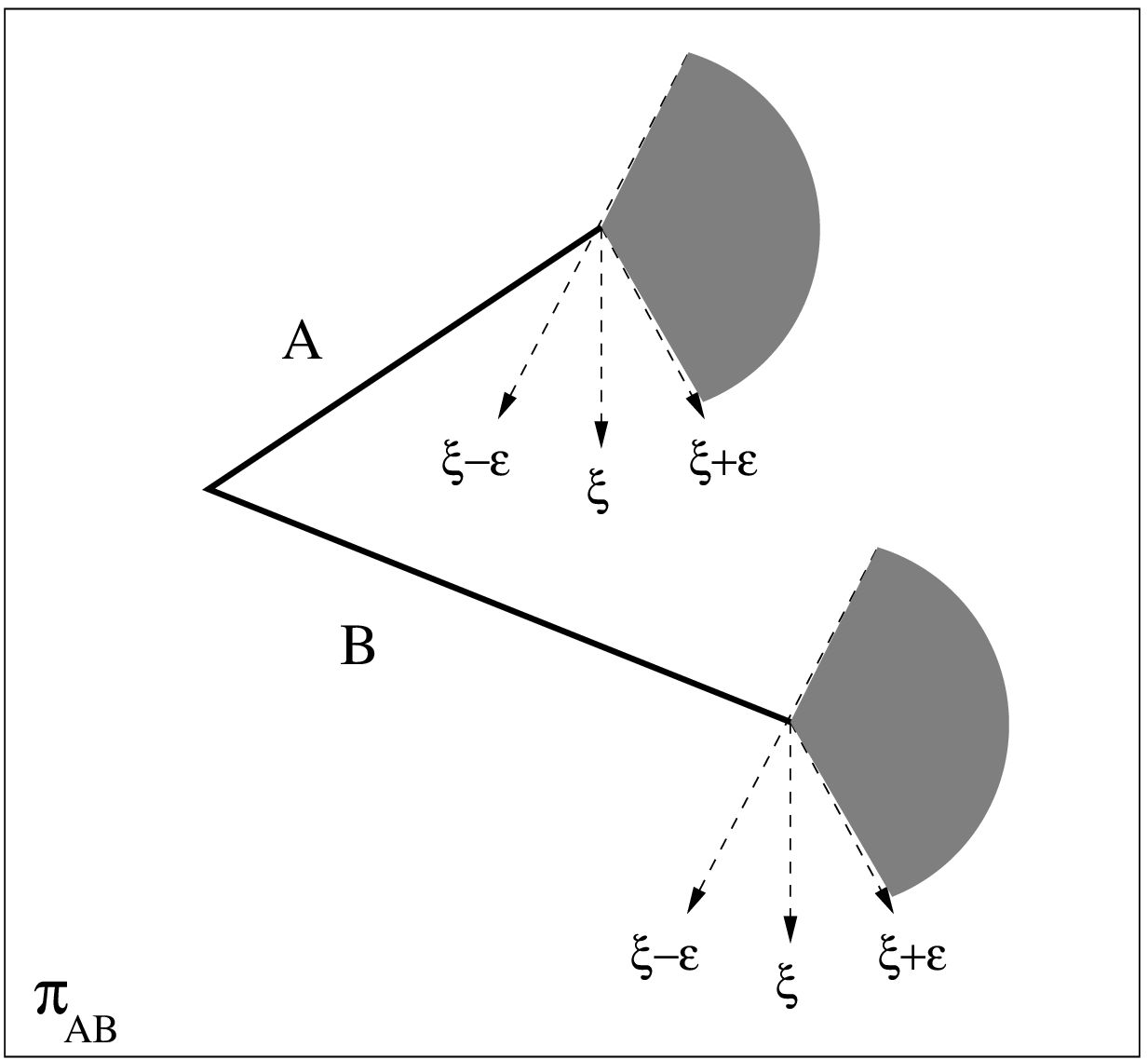,height=6cm}}
\vspace*{13pt}
\fcaption{}
\label{fig:cone}
\end{figure}

Let $N=\;]\xi - \epsilon ; \xi + \epsilon[\,$ be a neighborhood of $\xi$ in $\alpha$ such that the
indicatrix travels only once along $N$. 
One or several segments adjacent to $A$ or $B$ can remain in the plane 
$\pi_{AB}$ generated by $A$ and $B$, but these segments lie in the shaded area in figure
\ref{fig:cone}. Indeed, any segment outside the shaded area would make the indicatrix cover $N$ one
more time.

\begin{figure}[htbp]
\vspace*{13pt}
\centerline{\psfig{file=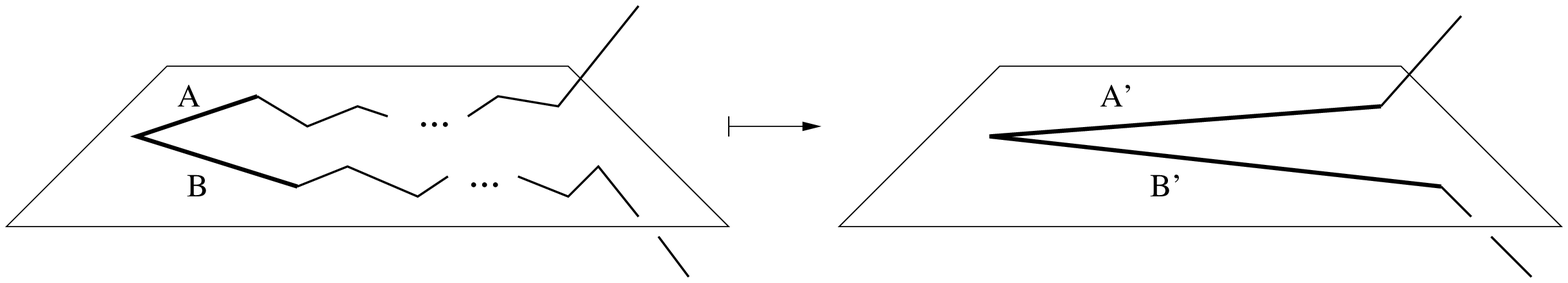,width=12cm}}
\vspace*{13pt}
\fcaption{}
\label{fig:trans}
\end{figure}

Clearly, the transformation illustrated in figure \ref{fig:trans} does not change $T(\xi_0)$ and 
$T(\xi_1)$; therefore, we only need to consider the case where the segments adjacent to $A$ and $B$ 
leave the plane $\pi_{AB}$. 

\begin{figure}[htbp] 
\vspace*{13pt}
\centerline{\psfig{file=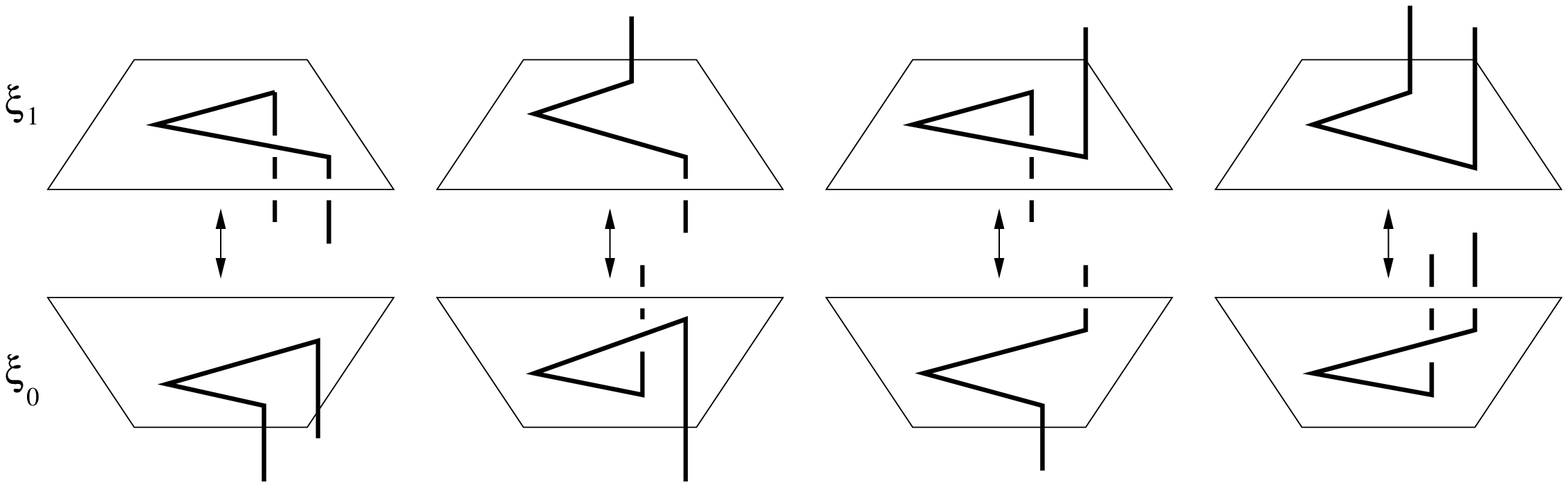,width=12cm}}
\vspace*{13pt}
\fcaption{}
\label{fig:ind}
\end{figure}

Let us check all the different possibilities. By proposition 1, the four cases 
illustrated in figure \ref{fig:ind} are sufficient. They correspond to the first Reidemeister move,
and we see that $T(\xi_1)=T(\xi_0)+1$.
\hfill \qed\kern0.8pt \medskip 

\noindent {\bf Corollary.} {\it Let $\xi,\xi\prime \in S^2\backslash I$; then $T_K(\xi)=
T_K(\xi\prime)+n$, where $n$ stands for the intersection number of the indicatrix with a path 
joining $\xi$ to $\xi\prime$.\hfill \qed\kern0.8pt}\medskip

The previous results were stated for $PL$ knots; nevertheless, they remain true for a wider 
class of knots.\medskip

\noindent {\bf Definition~6.} {\it A piecewise $C^2$ knot is the image of a closed space curve 
$\gamma \colon [0;1] \to {\R}^3$ twice continuously differentiable everywhere except on a finite 
number of points, satisfying for all $t_0 \in [0;1]$:
\begin{alphlist}
\item $\lim_{t \to t_0^{+}} \dot \gamma(t)$ and $\lim_{t \to t_0^{-}} \dot \gamma(t)$ exist and are 
non-zero;
\item $\lim_{t \to t_0^{+}} {\dot \gamma(t)\over ||\dot\gamma(t)||} \neq 
- \lim_{t \to t_0^{-}} {\dot\gamma(t) \over ||\dot\gamma(t)||}$.
\end{alphlist}
}\medskip

The indicatrix of a piecewise $C^2$ knot $K$ is defined in the obvious way, and it is always 
possible to approach $K$ with a sequence $\{K_n\}$ of $PL$ knots such that $I_{K_n}
\rightarrow I_K$. Hence, proposition 1 and 2 are true for piecewise $C^2$ knots, as well as
the corollary.\medskip

\begin{figure}[htbp] 
\vspace*{13pt}
\centerline{\epsfig{file=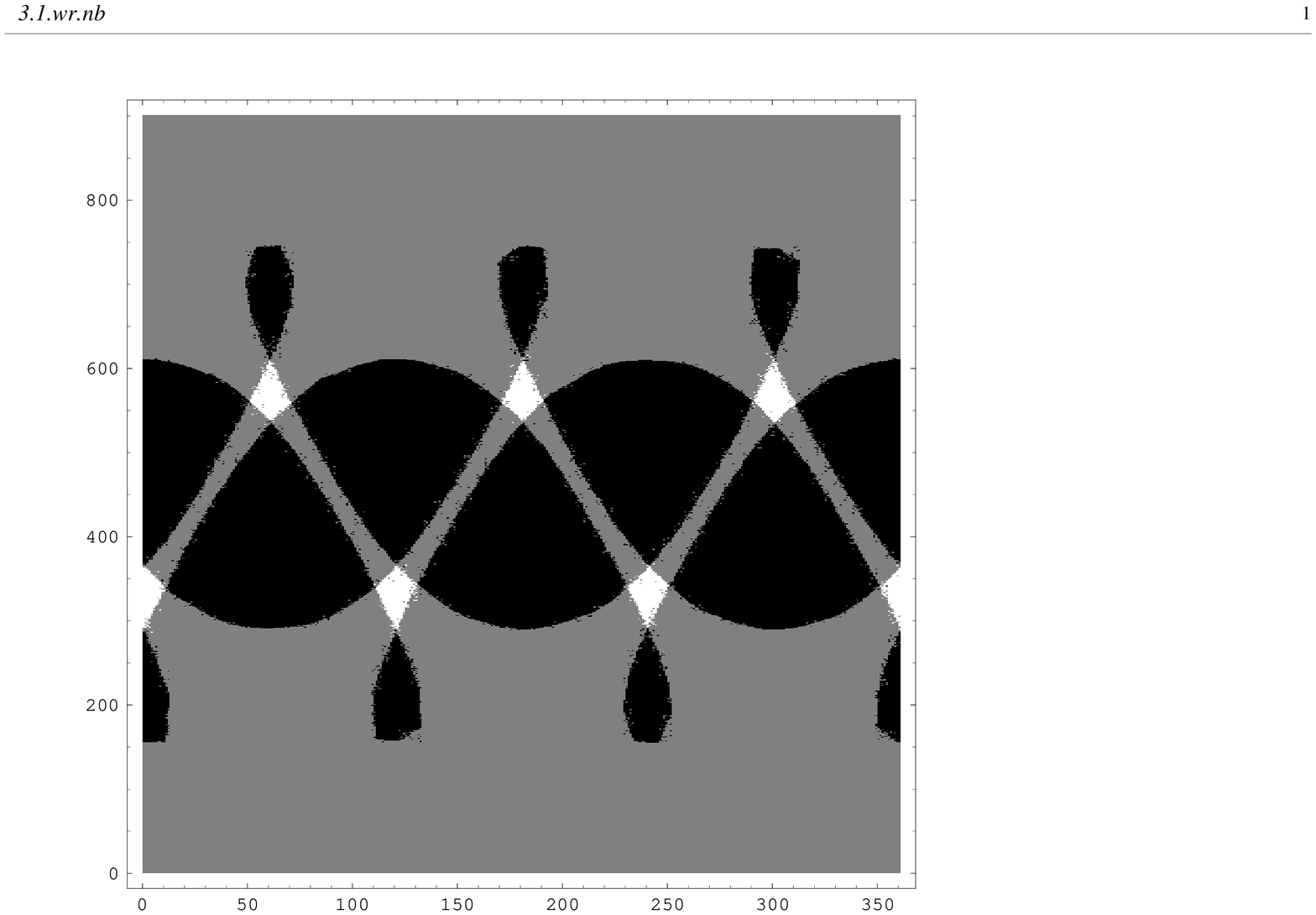,height=7cm,clip=,bbllx=70,bblly=90,bburx=510,bbury=540,
angle=-90}} 
\vspace*{13pt}
\fcaption{}
\label{fig:Akos}
\end{figure}

\noindent {\bf Illustration.} The picture shown on figure \ref{fig:Akos} was obtained by Akos Dobay
at the University of Lausanne. It represents $S^2$ via cylindrical coordinates with $360 \times 900$ 
evaluations of the function $T_K$ associated with a given smooth trefoil knot. For this particular 
space curve, these calculations provide a kind of ``experimental check'' of our results.  

\section{The Writhe of a Knot}
Let $K$ be a piecewise $C^2$ knot in ${\R}^3$; since $T_K(\xi)=T_K(-\xi)$ for all $\xi$ in 
$S^2\backslash I$,
$$ Wr(K)=\frac{1}{4\pi} \int\!\!\!\int_{S^2} T_K(\xi)\,d\xi = 
\frac{1}{2\pi} \int\!\!\!\int_{\frac{1}{2}S^2} T_K(\xi)\,d\xi.$$
By the corollary, we get the following formula, related to proposition 4 from Aldinger 
{\tenit et al.\/}:\medskip

\noindent {\bf Proposition~3.} {\it Let $K$ be a piecewise $C^2$ knot, $I$ its indicatrix,
$A_0,A_1,\dots,A_r$ the connected components of a hemisphere minus $I$, and $T_K(\xi_0)$ the
Tait number of $K$ relatively to some $\xi_0 \in A_0$. Then, the writhe of $K$ is given by
$$Wr(K)=T_K(\xi_0) + \frac{1}{2\pi}\sum_{i=1}^r n_i\cdot area(A_i) ,$$
where $n_i$ stands for the intersection number of the indicatrix with a path joining $A_0$ 
to $A_i$.\hfill \qed\kern0.8pt}\medskip

\begin{figure}[htbp]
\vspace*{13pt}
\centerline{\psfig{file=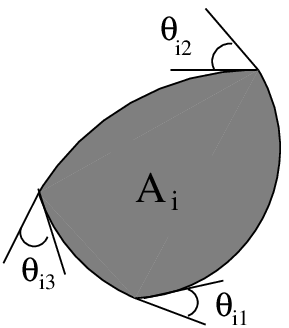}}
\vspace*{13pt}
\fcaption{}
\label{fig:angles}
\end{figure}

Let us now suppose that $K$ is a $PL$ knot. In this case, the $A_i$ are domains of $S^2$ 
delimited by arcs of great circles, that is, geodesics. By Gauss-Bonnet: $area(A_i)=2\pi-\sum_j 
\theta_{ij}$, where the $\theta_{ij}$ are the exterior angles of $A_i$ (see figure \ref{fig:angles}).
Given $K$, the computation of these angles is very easy. With the notations of the previous
proposition, we get:\medskip

\noindent {\bf Proposition~4.} {\it The writhe of a $PL$ knot $K$ is given by }
$$Wr(K)=T_K(\xi_0) + \frac{1}{2\pi}\sum_{i=1}^r n_i\cdot (2\pi-\sum_j  \theta_{ij}).$$
\vskip-1.1cm \hfill\qed\kern0.8pt
\vskip1.1cm 
For lattice knots, the calculation is immediate. \medskip
 
\noindent {\bf Definition~6.} {\it A lattice knot is a $PL$ knot on a cubic lattice in 
${\R}^3$.}\medskip

Let $K$ be a lattice knot; its indicatrix divides a hemisphere into four connected
components $A_1,\dots,A_4$, each of area $\frac{\pi}{2}$. Hence:
$$Wr(K)=\frac{1}{2\pi} \sum_{i=1}^4 T_K(A_i)\cdot \frac{\pi}{2} = \frac{T_K(A_1)+T_K(A_2)+
T_K(A_3)+T_K(A_4)}{4}.$$
Using only the first proposition, we have proved:\medskip

\noindent {\bf Proposition~5.} {\it If $K$ is a lattice knot, then $4 \cdot Wr(K)$ is an 
integer.\hfill \qed\kern0.8pt}\medskip

Furthermore, using the second proposition, it is an easy exercise to implement a program 
computing the writhe, given a diagram of the lattice knot.

\nonumsection{Final Remark}
The average crossing number of a $PL$ knot $K\/$ is defined by 
$$A(K) = \frac{1}{4\pi}\int\!\!\!\int_{S^2}|T_K(\xi)|\,d\xi.$$ 
Our method also applies to the computing of this number, leading to a closed formula. 
Here, all the curves on $S^2$ corresponding to the second Reidemeister move have to be taken 
into account, as well as the indicatrix $I$.

\nonumsection{Acknowledgements}
The author wishes to express his thanks to Claude Weber for his help, to Akos Dobay for 
his remarkable picture, and to Mathieu Baillif.

\nonumsection{References}

\end{document}